\documentclass[11pt]{amsart}
\usepackage{epsfig,psfrag,verbatim,amssymb,amsfonts}
\newtheorem{theorem}{Theorem}
\newtheorem{lemma}[theorem]{Lemma}
\newtheorem{prop}[theorem]{Proposition}

\newtheorem{theo}{Theorem}

\newtheorem*{rem}{Remarks}

\newcommand{\field}[1]{\mathbb{#1}}
\newcommand{\SL}{\mathcal{S}}
\newcommand{\R}{\field{R}}

\newcommand{\Z}{\field{Z}}
\newcommand{\N}{\field{N}}

\newcommand{\PP}{\field{P}}
\newcommand{\wPP}{\widetilde{\PP}}
\newcommand{\wP}{\widetilde{P}}
\newcommand{\bPP}{\bar{\PP}}
\newcommand{\bP}{\bar{P}}
\newcommand{\EE}{\field{E}}
\newcommand{\wEE}{\widetilde{\EE}}
\newcommand{\wE}{\widetilde{E}}

\newcommand{\om}{\omega}
\newcommand{\Om}{\Omega}
\newcommand{\eps}{\varepsilon}

\newcounter{constante}
\newcommand{\con}[1]{
\immediate\write 1{\noexpand\newlabel{#1}{{\theconstante}{\theconstante}}}
                    c_{\theconstante}
                    \stepcounter{constante}
                   }
\newcommand{\abel}[1]{}

\begin{document}

\setcounter{page}{1}

\title[Excited random walks]
{Recurrence and transience of excited random walks on $\Z^d$ and strips}

\maketitle

\vspace*{-2mm}
\begin{center}
\textit{Dedicated to the memory of Prof.\ Hans G.\ Kellerer (1934 - 2005)}\vspace{5mm}

{\sc By  Martin P.W.\ Zerner
}
\end{center}\vspace*{5mm}

{\footnotesize {\sc Abstract}. We investigate excited random walks on $\Z^d, d\ge 1,$ and on planar strips
$\Z\times\{0,1,\ldots,L-1\}$ which have a drift in a given direction. The strength of the drift may
depend on a random i.i.d.\ environment and on the local time of the walk.
We give exact criteria for recurrence and transience, thus generalizing
results by Benjamini and Wilson for once-excited random walk on $\Z^d$ and by the author for 
multi-excited random walk on $\Z$.
}\vspace*{5mm}

\section{Introduction}
We consider excited random walks (ERWs), precisely to be defined below, which move
on either $\Z^d$ or strips, i.e.\ which have state space
\[Y=\Z^d\quad (d\ge 1)\quad\mbox{or}\quad Y=\Z\times\{0,1,\ldots,L-1\}\subset \Z^2\quad (L\ge 2).\]
In general, ERWs are not Markovian. Instead, the transition probabilities may depend
on how often the walk has previously visited its present location and additionally
on the environment at this location.

To be more precise, let us first fix two quantities for the rest of the paper:
A direction $\ell$ and the so-called ellipticity constant $\kappa$.
In the 
case $Y=\Z$ or $Y=\Z\times\{0,1,\ldots,L-1\}$ we always choose 
$\ell=e_1\in Y$ to
be the first standard unit vector. In the case $Y=\Z^d$, $d\ge 2$, we let $\ell\in\R^d$ be any direction 
with $|\ell|_1=1$.
The ellipticity constant $\kappa\in(0, 1/(2d)]$ will be a uniform lower bound for the probability of the walk to jump from $x$ to any nearest neighbor of $x$.
Then an environment $\om$ for an ERW is an element of 
\begin{eqnarray*}
\Om&:=&\bigg\{\left(\left(\left(\om(x,e,i)\right)_{|e|=1}\right)_{i\geq 1}\right)_{x\in Y}
\in
   [\kappa,1-\kappa]^{2d\times\N\times Y } \bigg|\ \forall x\in Y \forall i\ge 1  \\ 
&&\hfill\sum_{e\in\Z^d,|e|=1}\om(x,e,i)=1,\
  \sum_{e\in\Z^d,|e|=1}\om(x,e,i)\, e\cdot \ell\ge 0\bigg\}.
\end{eqnarray*}
Here in the case of $Y$ being a strip, $d=2$ and $x+e$ is modulo $L$ in the second coordinate.

An ERW starting at $x\in Y$ in an environment $\om\in\Om$ is an $Y$-valued
process $(X_n)_{n\ge 0}$ on some suitable probability space $(\Om',\mathcal F, P_{x,\om})$
for which the history process $(H_n)_{n\ge 0}$ defined by $H_n:=(X_m)_{0\le m\le n}\in Y^{n+1}$
is a Markov chain which satisfies $P_{x,\om}$-a.s.\
\begin{eqnarray*}
P_{x,\om}[X_0=x]&=&1,\\
P_{x,\om}[X_{n+1}=X_n+e\mid H_n]&=& \om(X_n,e,\#\{m\le n\mid X_m=X_n\}).
\end{eqnarray*}
Thus $\om(x,e,i)$ is the probability to jump upon the $i$-th visit to $x$ from $x$ to $x+e$.
In the language introduced in \cite{ze}, an environment $\om\in\Om$ consists of infinite sequences 
of cookies attached to each site $x\in Y$. The $i$-th  cookie at $x$ is the transition 
vector $(\om(x,e,i))_{|e|=1}$ to the neighbors $x+e$ of $x$. Each time the walk visits $x$ it 
removes the first cookie from the sequence of cookies at $x$ and then jumps according to this cookie to a neighbor of $x$. Note that the assumption 
$\sum_{e}\om(x,e,i)e\cdot \ell\ge 0$ means that 
we allow only cookies which create a non-negative drift in direction $\ell$.
A model in which different sites may induce drift into opposite directions has been studied in \cite{abk05}.

The model described above generalizes ERW as introduced by  Benjamini and Wilson \cite{bewi}.
Their walk, which we will call BW-ERW, is an ERW on $\Z^d$, $d\ge 1$, in the environment $\om$ given by  
$\om(x,e,i)=1/(2d)$ for all $(x,e,i)$ with the only exception that
$\om(x,\pm e_1,1)=1/(2d)\pm\eps$, where $0<\eps<1/(2d)$ is fixed. Thus on the first visit to any site $x$, BW-ERW steps to $x\pm e_1$ with probability $1/(2d)\pm\eps$
and to all the other neighboring sites $x+e$ with probability
$1/(2d)$, while on any subsequent visit to $x$  a neighbor is chosen uniformly at random.
A main result of \cite{bewi} is the following.
\begin{theo}\label{d2}{\rm (see \cite{bewi})} 
BW-ERW on $\Z^d$, $d\ge 2$,
is transient in direction $e_1$, i.e.\ $X_n\cdot e_1\to\infty$
almost surely as $n\to\infty$.
\end{theo}\abel{d2}
Besides this it is also shown in \cite{bewi} that BW-ERW has positive liminf speed if $d\ge 4$. Kozma extended this result to $\Z^3$ in \cite{K03} and very recently even to $\Z^2$ in \cite{K05}. 

Unfortunately, we were not able to adapt the technique of proof introduced in \cite{bewi}, which was also used in \cite{K05}, to the more general
setting described above as we shall explain now.

Firstly, the proof in \cite{bewi}  relies on coupling the BW-ERW $(X_n)_n$
to a simple
symmetric random walk $(Y_n)_n$ in such a way that
$0\le (X_n-Y_n)\cdot e_1$ is non-decreasing in $n$ 
and $X_n\cdot e_i=
Y_n\cdot e_i$ for all $n$ and $i\ge 2$.
It is not clear to us how to achieve such a coupling
if one allows the drift 
to point into a direction
other than a coordinate direction $e_i$, e.g.\ by letting 
$\om(x,\pm e_j,1)=1/(2d)\pm\eps$ for $j=1,2$ and $\om(x,e,i)=1/(2d)$ for all
other $(x,e,i)$.

Secondly, the proof uses so-called tan points, which are points $x$ to the right
of which no other point has been visited prior to $x$. It is easy to see that
any tan point for $(Y_n)_n$ is also a tan point for $(X_n)_n$. Moreover, any
tan point for $(X_n)_n$ will be a point with a cookie, which will push $(X_n)_n$ to the right.
Then, roughly speaking,  using a lower bound on the number of tan points for
$(Y_n)_n$, one gets a lower bound on the number of cookies consumed by $(X_n)_n$, which
Benjamini and Wilson show to be sufficient to ensure transience to the right.
We do not see how this line of proof could be adapted to other settings, for instance if 
the excitement occurs not on the first but only on the second visit to a site.

In the present paper we suggest an alternative method of proof, based on martingales and on
the environment viewed from the particle, which applies to BW-ERW and  to other more 
general settings, in which $\om$ is sampled from $\Om$ according to a probability measure 
$\PP$ on $\Om$
 such that the family
\begin{equation}\label{scr}
(\om(x,\cdot,\cdot))_{x\in Y}\quad\mbox{is i.i.d.\ under $\PP$}.
\end{equation}
Throughout the paper we will assume (\ref{scr}) and denote the expectation with respect to $\PP$ by
$\EE$.
Note that we do not assume independence between different cookies at the same site nor between transition probabilities to different neighbors of the same site, but only between cookies at different sites. 
An important quantity will be the total drift $\delta^x$ in direction 
$\ell$ of all the cookies
stored at  site $x\in Y$, i.e.
\[\delta^x(\om):=\sum_{i\ge 1, |e|=1}\om(x,e,i)\, e\cdot \ell\, .\]
Note that by definition of $\Om$, $\delta^x(\om)\ge 0$ for all $x\in Y$ and  $\om\in\Om$.
We shall generalize Theorem \ref{d2} as follows.
\begin{theorem} \label{dd} Let $d\ge 2$, $Y=\Z^d$
and $\EE[\delta^0]>0$. 
Then the walk is for $\PP$-almost all $\om$ transient in direction $\ell$, i.e.\
 $P_{0,\om}$-a.s.\ $X_n\cdot \ell\to \infty$ as $n\to\infty$.
\end{theorem}\abel{dd}
The technique of proof improves methods  used  in \cite{ze} to show the following result for $d=1$. Some simulation studies for $Y=\Z$ can be found in \cite{anre}.
\begin{theo}\label{d1}{\rm (see \cite[Theorem 12]{ze})} Let $Y=\Z$.
Then for $\PP$-almost all environments $\om\in\Om$, $(X_n)_n$ is recurrent, 
i.e.\ returns $P_{0,\om}$-a.s.\ (infinitely often)
to its starting point, if and only if $\EE[\delta^0]\le 1$.
\end{theo}\abel{d1}
In fact, \cite[Theorem 12]{ze} is more general since it does not need 
any ellipticity condition and allows the environment to be stationary and ergodic only instead of i.i.d..
In the present paper we shall generalize Theorem \ref{d1} to strips as follows.
\begin{theorem}\label{kellerer} Let $Y=\Z$ and $L=1$ or 
$Y=\Z\times\{0,\ldots,L-1\}$ for some $L\ge 2$.
If $\EE[\delta^0]>1/L$ then the walk is  for $\PP$-almost all $\om$ transient in direction $e_1$.
If $\EE[\delta^0]\le 1/L$ then the walk is for $\PP$-almost all $\om$ recurrent,
and moreover $P_{0,\om}$-a.s.\ 
 $\limsup_{n\to\infty}X_n\cdot e_1=\infty$ and $\liminf_{n\to\infty}X_n\cdot e_1=-\infty$.
\end{theorem}\abel{kellerer}
So if the strip is made wider and wider while the distribution of 
$\om(x,\cdot,\cdot)$ is kept fixed, the walk will eventually become 
transient if $\EE[\delta^0]>0$. This provides some additional support for Theorem \ref{dd}.
\section{Preliminaries}
For $z\in\R,n\in\N\cup\{\infty\}$ we let
\[ D_n^z:=\sum_{x\in \SL_z}
\sum_{i=1}^{\#\{m<n\ \mid\  X_m=x\}}\sum_{|e|=1}\om(x,e,i) e\cdot \ell \]
denote the drift absorbed by the walk by time $n$ while visiting the slab $\SL_z:=\{x\in Y\mid z\le x\cdot\ell<z+1\}$. 
Then $D_n:=\sum_{z\in\Z}D_n^z$ is the total drift encountered by the walk up to time $n$.
Observe that $D_n^z\ge 0$ and therefore also $D_n\ge 0$ for all $\om\in\Om$ and all paths $(X_m)_m$.

By standard arguments, for any $\om\in\Om$ the process $(M_n)_{n\ge 0}$ defined by
\begin{equation}\label{meyer}
M_n:=X_n\cdot \ell-D_n
\end{equation}\abel{meyer}
is a  martingale under $P_{0,\om}$ with respect to the filtration
generated by $(X_n)_{n\ge 0}$. Indeed, (\ref{meyer}) is just the Doob-Meyer decomposition
of the submartingale $(X_n\cdot\ell)_n$.

In the setting considered in \cite{bewi} and \cite{ze} part of the following fact was achieved by coupling
 the  ERW to a simple symmetric random walk staying always 
to the left of the ERW. For the present more general setting we need a different argument.
\begin{lemma}\label{beata}
Let $\om\in\Omega$. Then $P_{0,\om}$-a.s.\ 
\[\liminf_{n\to\infty} X_n\cdot\ell\in\{-\infty,+\infty\}\quad\mbox{and}\quad
\limsup_{n\to\infty} X_n\cdot\ell=+\infty.\]
In particular, for all $x\ge 0$,
\[T_x:=\inf\{n\ge 0\mid X_n\cdot\ell\ge x\}<\infty\quad\mbox{$P_{0,\om}$-a.s..}\]
\end{lemma}\abel{beata}
\begin{proof}
It follows from ellipticity and the Borel-Cantelli lemma that $\liminf_{n} X_n\cdot\ell\notin\R$, which implies the first statement.

For the statement about $\limsup$, let $x\ge 0$.
Since $D_n\ge 0$ for all $n$, the martingale  $(M_{n\wedge T_x})_n$
is bounded from above by $x$ and hence converges $P_{0,\om}$-a.s.\ to a finite limit
as $n\to\infty$. 
Therefore, it suffices to show, that  $(M_n)_n$ itself $P_{0,\om}$-a.s.\ does not converge, because
then the convergence of  $(M_{n\wedge T_x})_n$ can only be due to $T_x$ being   $P_{0,\om}$-a.s.\ finite.

So if $(M_n)_n$ did converge, then
$|(X_{n+1}-X_n)\cdot\ell -(D_{n+1}-D_n)|\to 0$ as $n\to\infty$. However, this is impossible. 
Indeed, let $e_0\in Y$ be a unit vector which maximizes 
$e_0\cdot\ell$. Then due 
to ellipticity and the Borel-Cantelli lemma,
$|(X_{n+1}-X_n)\cdot\ell|=e_0\cdot\ell$ infinitely often, whereas, again by ellipticity, for all $n$ 
and some random $i=i(n)\in\N$,
\begin{eqnarray*}
|D_{n+1}-D_n|
&=&\bigg|\sum_e\om(X_n,e,i)e\cdot\ell\bigg|\\
&\le& \bigg(|\om(X_n,e_0,i)-\om(X_n,-e_0,i)|+\sum_{e\ne\pm e_0}\om(X_n,e,i)\bigg)\ e_0\cdot \ell\\
&=&\left(1-2\left(\om(X_n,e_0,i)\wedge \om(X_n,-e_0,i)\right)\right)\ e_0\cdot\ell,
\end{eqnarray*}
which is at most $(1-2\kappa)e_0\cdot \ell$.
\end{proof}
\begin{lemma}\label{beatus}
For all $\om\in\Omega$ and all $x\ge 0$,
$E_{0,\om}[D_{T_x}]\leq x+1$.
\end{lemma}\abel{beatus}
\begin{proof}
By the Optional Stopping Theorem for all $n\in\N$,
$0=E_{0,\om}[M_{T_x\wedge n}]$ and consequently by (\ref{meyer}), 
$E_{0,\om}[D_{T_x\wedge n}]= E_{0,\om}[X_{T_x\wedge n}\cdot \ell]\le x+\max_e e\cdot\ell\le x+1.$
The statement now follows from monotone convergence.
\end{proof}

Now we introduce some notation taken from \cite{ze} for the 
cookie environment left over by the walk.
For $\om\in\Om$ and any finite sequence 
$(x_n)_{n\le m}$ in $Y$ we define $\psi(\om,(x_n)_{n\leq m})\in\Om$ by
\[
\psi(\om,(x_n)_{n\leq m})(x,e,i):=\om\left(x,e,i+\#\{n<m\mid x_n=x\}\right).
\]
This is the environment created by the ERW by following the path $(x_n)_{n\le m}$ and removing all the first cookies encountered, 
except for the last visit to $x_m$.
Finiteness of $T_1$, guaranteed by Lemma \ref{beata}, implies that the Markov transition kernel 
\[R(\om,\om'):=P_{0,\om}\left[\theta^{X_{T_1}}\left(\psi\left(\om, H_{T_1}\right)\right)=\om'\right]\]
for $\om,\om'\in\Om$ is well-defined. Here $\theta^z$ denotes the spatial shift of the 
environment by $z$, i.e.\ $\theta^z(\om(x,\cdot,\cdot)):=\om(x+z,\cdot,\cdot)$. The probability measure $R(\om,\cdot)$ is
the distribution of the modified environment $\om$
viewed from the particle at time $T_1$. Note that it
is supported on those countably many $\om'\in\Omega_\ell$, which are obtained from $\om$ by 
removing finitely many cookies from $\om$.
\begin{lemma}\label{feller}
$R$ is weak Feller, i.e.\ convergence w.r.t.\ the product topology on $\Om$ of $\om_n\in\Om$
towards $\om\in\Omega$ as $n\to\infty$ implies
\begin{equation}\label{sums}
\left|\sum_{\om'\in\Om}R(\om_n,\om')f(\om')- \sum_{\om'\in\Om}R(\om,\om')f(\om')\right|
\longrightarrow 0\quad\mbox{as $n\to\infty$}
\end{equation}\abel{sums}
for any bounded continuous function $f:\Om\to\R$.
\end{lemma}\abel{feller}
Note that due to the discreteness of $R(\om,\cdot)$ only countably many terms in the  sums in (\ref{sums})
do not vanish.
\begin{proof}
Let $\eps>0$. Since $T_1$ is $P_{0,\om}$-a.s.\ finite due to Lemma
\ref{beata}, there is some finite $t$  such that 
\begin{equation}\label{bat}
\eps\ >\ P_{0,\om}[T_1>t]\ =\ 1-\sum_{\pi\in\Pi_t}P_{0,\om}\left[\mbox{$(X_m)_m$ follows $\pi$}\right], 
\end{equation}\abel{bat}
where $\Pi_t$ denotes the set of nearest-neighbor paths $\pi$ starting at the origin and 
ending at time $T_1(\pi)$ with $T_1\le t$.
Since $\om_n\to\om$,
\begin{equation}\label{wahe}
P_{0,\om_n}\left[\mbox{$(X_m)_m$ follows $\pi$}\right]\longrightarrow 
P_{0,\om}\left[\mbox{$(X_m)_m$ follows $\pi$}\right]\qquad\mbox{as $n\to\infty$}
\end{equation}\abel{wahe}
for all $\pi\in\Pi_t$. Therefore, by (\ref{bat}), 
\begin{equation}\label{ura}
P_{0,\om_n}[T_1>t]<\eps\quad\mbox{for $n$ large}.
\end{equation}\abel{ura}
Now  partition $\Pi_t$ into sets $\Pi_t^z$ according to the final point $z$ of the paths.
Then the left-hand side of (\ref{sums}) can be bounded from above by
\begin{eqnarray}\nonumber
&&\sum_{z\in Y}\sum_{\pi\in\Pi_t^z}\bigg|P_{0,\om_n}[\mbox{$(X_m)_m$ follows $\pi$}] 
f\left(\theta^z(\psi(\om_n,\pi))\right)\\
&&\hspace*{19mm}-\ P_{0,\om}[\mbox{$(X_m)_m$ follows $\pi$}] 
f\left(\theta^z(\psi(\om,\pi))\right)\bigg|\label{edi}\\ \nonumber
&&+\ c P_{0,\om}[T_1>t]+ c P_{0,\om_n}[T_1>t],
\end{eqnarray}\abel{edi}
where $c$ is a bound on $|f|$.
Since $f$ is continuous, $f\left(\theta^z(\psi(\om_n,\pi))\right)$ converges to
$f\left(\theta^z(\psi(\om,\pi))\right)$ as $n\to\infty$. 
Together with (\ref{bat}), (\ref{wahe}) and  (\ref{ura}) this
shows that the whole expression in (\ref{edi}) is less than $2c\eps$ for $n$ large. 
\end{proof}
\begin{lemma}\label{compact}
There is a probability measure $\wPP$ on $\Om$ 
which is invariant under
$R$ and under which 
\begin{equation}\label{steiff}
(\om(x,\cdot,\cdot))_{x\in Y, x\cdot\ell\ge 0}  \quad\mbox{has the same distribution 
as under $\PP$.}
\end{equation}
\end{lemma}\abel{compact}
\begin{proof}
Being  a closed subset of the compact set $[\kappa,1-\kappa]^{2d\times\N\times Y}$,
 $\Om$ is compact, too.
Consequently, the set of all probability measures on $\Om$ is compact as well.
Since the set $\mathcal M$ of all probability measures on $\Om$ 
under which (\ref{steiff}) holds
is a closed subset of this compact set, $\mathcal M$ is compact, too. 
Moreover, observe that $MR\in\mathcal M$ for  all $M\in\mathcal M$ since the part of the
environment $\psi(\om, H_{T_1})$ which is to the right
of  $X_{T_1}$ has by time $T_1$ not been touched by the walk yet and is therefore still i.i.d.. 
Hence, since $R$ is weak Feller due to Lemma \ref{feller} the 
statement follows from standard arguments, see e.g.\ \cite[Theorem 12.0.1 (i)]{metw}.
\end{proof}
For the remainder of this paper we fix $\wPP$ according to Lemma \ref{compact} and
let $\wEE$ be its expectation operator. We also introduce the annealed probability measures 
$P_0=\PP\times P_{0,\om}$ and 
$\wP_0=\wPP\times P_{0,\om}$ with expectation operators $E_0$ and $\wE_0$, respectively, 
which one gets by averaging the so-called quenched measure $P_{0,\om}$ over $\EE$ and 
$\wEE$, respectively, i.e.\ $P_0[\cdot]=\EE[P_{0,\om}[\cdot]]$ and 
$\wP_0[\cdot]=\wEE[P_{0,\om}[\cdot]]$. 
The following statement is similar to \cite[Lemma 11]{ze}.
\begin{lemma}\label{blanda}
If $Y$ is a strip or $\Z$ then 
$\wE_0[D_\infty^0]\le 1$.
If $Y=\Z^d,\ d\ge 2,$ then
$\wE_0[D_\infty^0]\le 2$.
\end{lemma}\abel{blanda}
\begin{proof}
Consider the stopping times defined by $\tau_0:=0$ and $\tau_{n+1}:=\inf\{n>\tau_n\ :\ 
X_n\cdot \ell\ge X_{\tau_n}\cdot \ell+1\}$ for $n\ge 0$.
Note that 
\begin{equation}\label{sch}
\tau_n=T_n\quad\mbox{if $Y$ is a strip or $\Z$ and}\quad \tau_n\le T_{2n}\quad\mbox{if $Y=\Z^d,\ d\ge 2,$}
\end{equation}\abel{sch}
because in the second case, due to $|\ell|=1$, $X_{\tau_{n+1}}\cdot \ell\le X_{\tau_{n}}\cdot \ell+2$. Since the slabs $\SL_{X_{\tau_n}\cdot\ell}$,\ $n\ge 0,$ are disjoint, we have $D_{T_K}\ge \sum_{n\ge 0} D_{T_K}^{X_{\tau_n}\cdot\ell}$ for all $K\ge 0$.
Therefore,
for all $0\le k<K/2$,
\begin{equation}\label{prince}
D_{T_K}\ge  \sum_{n=0}^{m}
 D_{\tau_{n+k}}^{X_{\tau_n}\cdot \ell},
\end{equation}\abel{prince}
where $m=m(K,k):=K-k$ for $Y$ being a strip or $\Z$ and $m(K,k):=\lfloor K/2\rfloor-k$ for $Y=\Z^d,\ d\ge 2$.
Indeed, in both cases $\tau_{n+k}\le T_K$ for
all $n\le m$ due to (\ref{sch}).
Consequently, by Lemma \ref{beatus} and (\ref{prince}),
\begin{equation}\label{paul}
K+1\ge \wE_0[D_{T_K}]\ge  \sum_{n=0}^{m}
 \wE_0\left[D_{\tau_{n+k}}^{X_{\tau_n}\cdot \ell}\right].
\end{equation}\abel{paul}
By conditioning on the history up to time $\tau_n$ and using the strong Markov property we get
\begin{eqnarray}
 \wE_0\left[D_{\tau_{n+k}}^{X_{\tau_n}\cdot \ell}\right]
& =& \nonumber
\wEE\left[
E_{0,\om}\left[E_{0,\theta^{X_{\tau_n}}\left(\psi(\om,H_{\tau_n})\right)}[ D_{\tau_{k}}^0]\right]\right]\\
&=&\ \label{nonumber}
 \wEE\left[\sum_{\om'\in\Om}E_{0,\om}\left[ E_{0,\om'}[D_{\tau_{k}}^0],
\theta^{X_{\tau_n}}\left(\psi(\om,H_{\tau_n})\right)=\om'\right]
\right]\\
&=&\ \wEE\left[\sum_{\om'\in\Om}E_{0,\om'}[ D_{\tau_{k}}^0] R^n(\om,\om')
\right]\nonumber 
\ =\  \wEE\left[E_{0,\om}[ D_{\tau_{k}}^0]\right],
\end{eqnarray}\abel{nonumber}
where $R^n$ denotes the $n$-th iteration of $R$ and 
the last identity holds due to $\wPP R^n=\wPP$. Consequently, we obtain from (\ref{paul}) that
$\wE_0[ D_{\tau_{k}}^0]\le (K+1)/m(K,k)$. 
Letting $K\to\infty$ gives, for all $k\ge 0$,
$\wE_0[ D_{\tau_{k}}^0]\le 1$ for the strip and $\Z$ and 
$\wE_0[ D_{\tau_{k}}^0]\le 2$  for $\Z^d$, $d\ge 2$. 
Monotone convergence as $k\to\infty$
then yields the claim.
\end{proof}
\section{Transience on $\Z^d$ and strips}
We denote by
\[A_\ell\ :=\ \left\{\lim_{n\to\infty}X_n\cdot\ell=+\infty\right\}\quad\mbox{and}\quad
B_\ell\ :=\ \left\{\forall n\ge 1\ X_n\cdot\ell>X_0\cdot \ell\right\}\]
the event that the walk tends to the right and the event that it 
stays forever strictly to the right of its initial point, respectively.
As a preliminary result, we are now going to prove 
Theorem \ref{dd} with $\wPP$ instead of $\PP$.
\begin{lemma}\label{finite}
Let $d\ge 2$,\ $Y=\Z^d$ and 
$\EE[\delta^0]>0$. Then 
$\wP_{0}[A_\ell]=1$.
\end{lemma}\abel{finite}
\begin{proof}
On $A_\ell^c$ the walk changes 
 sign $\wP_0$-a.s.\ infinitely often
 due to Lemma \ref{beata}.
 Therefore, because of ellipticity, on $A_\ell^c$ it
 also visits $\wP_0$-a.s.\ infinitely many sites in the slab 
$\SL_0$.
Among these sites $x$ there are $\wP_0$-a.s.\ infinitely many ones with 
$\sum_{|e|=1,i\le I}\om(x,e,i)e\cdot\ell>\eps$ for some
$\eps>0$ and some finite $I$ due to  the assumption of independence in the environment and $\EE[\delta^0]>0$. 
Again by ellipticity, on $A_\ell^c$, $\wP_0$-a.s.\ infinitely many of those sites will be visited at least $I$ times. This 
yields  that on $A_\ell^c$, $\wP_0$-a.s.\ 
$D_\infty^0=\infty$, which would contradict Lemma
\ref{blanda} unless $\wP_0[A_\ell^c]=0$.
\end{proof}
The following type of result is  standard, see e.g.\ \cite[Lemma 1]{S94}, \cite[Proposition 1.2]{SZ99}
 and \cite[Lemma 8]{ze}.
\begin{lemma}\label{D}
Let $\om\in\Omega$ such that $P_{0,\om}[A_\ell]>0$. Then   $P_{0,\om}[A_\ell\cap B_\ell]>0$.
\end{lemma}\abel{D}

\begin{proof}
By assumption there is a finite nearest-neighbor path $\pi_1$ starting
at $0$ and ending at some $a$ with $a\cdot\ell> d$ 
such that with positive
$P_{0,\om}$-probability the walk first follows  $\pi_1$ and then stays  
to the right of $a$, while tending to the right, i.e.
\begin{equation}\label{forch}
P_{0,\om}[\mbox{$(X_n)_n$ follows $\pi_1$}]\ P_{a,\psi(\om,\pi_1)}[A_\ell\cap B_\ell]>0,
\end{equation}\abel{forch}
see Figure \ref{paths}.
\begin{figure}[t]
\epsfig{file=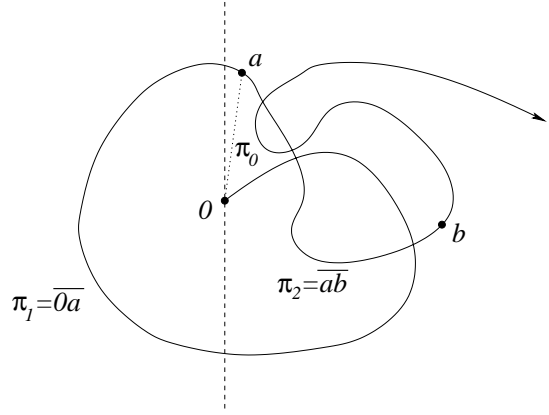,height=5.5cm, angle=0}
\caption{For the proof of Lemma \ref{D}. The path $\pi_1$ from 0 to $a$ is cut out and replaced
by the dotted path $\pi_0$.}
\label{paths}
\end{figure}\abel{paths}
In particular, the second factor in (\ref{forch}) is positive.
Now on $A_\ell$, the walk can visit sites on the path $\pi_1$ only finitely often. Therefore,
there is another path $\pi_2$ of length $m_2$ entirely to the right of $a$ which starts at $a$ and ends at some $b$ 
such that
\begin{eqnarray*}
0&<&P_{a,\psi(\om,\pi_1)}[\{\mbox{$(X_n)_n$ follows $\pi_2$}\}\cap \{\forall n\ge m_2\ X_n\notin\pi_1\}\cap A_\ell\cap B_\ell]\\
&\le&
P_{b,\psi(\om,(\pi_1,\pi_2))}[\{\forall n> 0\ X_n\cdot\ell> a\cdot\ell,  X_n\notin\pi_1\}\cap A_\ell].
\end{eqnarray*}
However,  on the event that the walk never visits $\pi_1$ the walk does not feel whether it moves in the environment
$\psi(\om,(\pi_1,\pi_2))$ or $\psi(\om,\pi_2)$.
Therefore, 
\[0<P_{b,\psi(\om,\pi_2)}[\{\forall n> 0\ X_n\cdot\ell> a\cdot\ell\}\cap A_\ell].
\]
Since $P_{a,\om}[\mbox{$(X_n)_n$ follows $\pi_2$}]>0$ due to ellipticity,
we get from this
\begin{eqnarray}\nonumber
0&<&P_{a,\om}[\mbox{$(X_n)_n$ follows $\pi_2$}]\  
P_{b,\psi(\om,\pi_2)}[\{\forall n> 0\ X_n\cdot\ell> a\cdot\ell\}\cap A_\ell]\\
&=&P_{a,\om}[\{\mbox{$(X_n)_n$ follows $\pi_2$}\}\cap A_\ell\cap B_\ell]
\ \le\ P_{a,\om}[A_\ell\cap B_\ell].\label{dumble}
\end{eqnarray}\abel{dumble}
Now because of $a\cdot\ell>d$ there is a nearest-neighbor path $\pi_0$
from $0$ to $a$ with $0<x\cdot\ell<a\cdot\ell$
for all sites $x$ on $\pi_0$ except for its starting and its end point.
By ellipticity, the walk will follow $\pi_0$ with positive $P_{0,\om}$-probability.
Therefore, due to (\ref{dumble}) and since 
$P_{a,\om}[A_\ell\cap B_\ell]=P_{a,\psi(\om,\pi_0)}[A_\ell\cap B_\ell]$,
\begin{eqnarray*}
0&<&P_{0,\om}[\mbox{$(X_n)_n$ follows $\pi_0$}]\ P_{a,\om}[A_\ell\cap B_\ell]\\
&=&
P_{0,\om}[\{\mbox{$(X_n)_n$ follows $\pi_0$}\}\cap A_\ell\cap B_\ell]
\end{eqnarray*}
by the strong Markov property. Hence $P_{0,\om}[A_\ell\cap B_\ell]>0$.
\end{proof}
We are now ready to prove a 0-1-law. We shall apply this result to $\bPP\in\{\PP,\wPP\}$.
\begin{prop}\label{01}
Let $\bPP$ be a probability measure on $\Om$ and
let $(\om(x,\cdot,\cdot))_{x\cdot\ell\ge 0}$  be i.i.d.\ under $\bPP$ . Then
$(\bPP\times P_{0,\om})[A_\ell]\in\{0,1\}$.
\end{prop}\abel{01}
\begin{proof} For short set $\bP_0=\bPP\times P_{0,\om}$.
Let us assume $\bP_0[A_\ell]>0$. We need to show $\bP_0[A_\ell]=1$. By Lemma \ref{D},
$\bP_0[B_\ell]>0$.
The following argument
is well-known, see e.g.\  \cite[Lemma 1.1]{SZ99} and \cite[Proposition 3]{ZM01}.
Fix $M\in \N$.
We define recursively  possibly infinite stopping times $(S_k)_{k\geq 0}$ and $(R_k)_{k\ge 0}$ by
$S_0:=T_M$,
\begin{eqnarray*}
R_k&:=& \inf\{n\ge S_k\mid X_n\cdot\ell<M\}\quad\mbox{and}\\
S_{k+1}&:=&\inf\left\{n\ge R_k\mid X_n\cdot\ell>\max_{m<n}X_m\cdot \ell\right\}.
\end{eqnarray*}
 Due to Lemma \ref{beata}
$S_0$ is $\bP_0$-a.s.\ finite and any subsequent 
$S_{k+1}$ is $\bP_0$-a.s.\ finite as well provided $R_k$ is finite.
Moreover, at each finite time $S_k$ the walk has reached a half space it has never 
touched before.
The environment $(\om(x+X_{S_k},\cdot,\cdot))_{x\cdot\ell\ge 0}$ in this half space 
is independent of the environment visited so far and has the same distribution as 
$(\om(x,\cdot,\cdot))_{x\cdot\ell\ge 0}$.
Hence the walk has probability $\bP_0[B_\ell]$ never to leave this half space again. Therefore,
by induction, $\bP_0[R_k<\infty]\le \bP_0[B_\ell^c]^k$, which goes to 0 as $k\to\infty$.
Consequently, there is a random integer $K$ with $R_K=\infty$. This means that
$X_n\cdot \ell\ge M$ for all $n\ge S_K$. 
Since this holds for all $M$, $\bP_0[A_\ell]=1$. 
\end{proof}
The following is the counterpart of Lemma \ref{finite} for $\Z$ and strips.
\begin{lemma}\label{finite2}
Let $Y=\Z$ and $L=1$ or $Y=\Z\times\{0,\ldots,L-1\}$ for some $L\ge 2$ and let $\EE[\delta^0]>1/L$. Then $\wP_0[A_\ell]=1$.
\end{lemma}\abel{finite2}
\begin{proof} Assume that $\wP_0[A_\ell]<1$. Then by Proposition \ref{01}, $\wP_0[A_\ell]=0$.
Therefore, the walk changes 
 sign $\wP_0$-a.s.\ infinitely often
 due to Lemma \ref{beata}.
However, if the walk crosses 
the finite set $\SL_0$ infinitely often then by ellipticity it
will eventually eat all the cookies in $\SL_0$, i.e.\ $\wP_0$-a.s.\ $D_\infty^0=\sum_{x\in \SL_0}\delta^x$. Hence
$\wE_0[D_\infty^0]=L\wE_0[\delta^0]>1$, which  contradicts
Lemma \ref{blanda}. 
\end{proof}
\begin{proof}[Proof of Theorem \ref{dd} and of transience in Theorem \ref{kellerer}]
By Lemma \ref{finite} and Lemma \ref{finite2}, respectively, 
$\wP_0[A_\ell]=1$. Therefore, due
to Lemma \ref{D},   $\wP_0[A_\ell\cap B_\ell]>0$. However, since (\ref{steiff}) holds under $\wPP$,
$P_0[A_\ell\cap B_\ell]=
\wP_0[A_\ell\cap B_\ell]>0$. Consequently, by Proposition \ref{01}, $P_0[A_\ell]=1$. 
\end{proof}

\section{Recurrence on strips}
\begin{proof}[Proof of recurrence in Theorem \ref{kellerer}]
Let $L\EE[\delta^0]\le 1$. We need to show that $P_0$-a.s.\ 
$\liminf_{n} X_n\cdot e_1\le 0$, since then, by ellipticity, $X_n=0$ infinitely
often. 
Assume the contrary.
 Then by Lemma \ref{beata}, $P_0[A_\ell]>0$. Consequently, by Lemma \ref{D},
even 
$P_0[A_\ell\cap B_\ell]>0$.
However, $P_0[A_\ell\cap B_\ell]=\wP_0[A_\ell\cap B_\ell]$.
Hence, by Proposition \ref{01}, 
\begin{equation}\label{unopoo}
\wP_0[A_\ell]=1.
\end{equation}\abel{unopoo} 
Now let $T_{-i}:=\inf\{n\mid X_n\cdot\ell\le -i\}$ for $i>0$. Then we have 
by the Optional Stopping Theorem for all $i,k,n\in\N$ and all $\om\in\Om$,
\begin{eqnarray}
0\ =\ E_{0,\om}[M_{T_k\wedge T_{-i}\wedge n}]\nonumber 
&=&k P_{0,\om}[T_k<T_{-i}\wedge n]-i P_{0,\om}[T_{-i}<T_k\wedge n]\nonumber \\
&&\label{domi}
+\ E_{0,\om}[X_n\cdot \ell, n<T_k\wedge T_{-i}]-E_{0,\om}[D_{T_k\wedge T_{-i}\wedge n}].
\end{eqnarray}\abel{domi} 
Using dominated convergence as $n\to\infty$ for both terms in (\ref{domi}), we obtain
\begin{eqnarray*}\frac{1}{k}E_{0,\om}[D_{T_k\wedge T_{-i}}]&=&
P_{0,\om}[T_k<T_{-i}]-\frac{i}{k}
P_{0,\om}[T_{-i}<T_k].
\end{eqnarray*}
Hence, due to (\ref{unopoo}), $\wPP$-a.s.\
$\lim_{i\to\infty}\lim_{k\to\infty}k^{-1}E_{0,\om}[D_{T_k\wedge T_{-i}}]=1$.
Splitting $D_n$ into $D_n^+:=\sum_{k\ge 0} D_n^k$ and $D_n^-:=\sum_{k< 0} D_n^k$
then yields
\begin{equation}\label{snape}
\lim_{i\to\infty}\lim_{k\to\infty}\frac{1}{k}E_{0,\om}[D^+_{T_k\wedge T_{-i}}]
= 1,
\end{equation}\abel{snape}
since $E_{0,\om}[D^-_{T_k\wedge T_{-i}}]\le \sum_{-i< x\cdot e_1< 0}\delta^x(\om)$,
which is $\wPP$-a.s.\ finite, does not depend on $k$ and thus vanishes when divided by $k\to\infty$.
However,
\begin{equation}\label{harry}
E_{0,\om}[D^+_{T_k\wedge T_{-i}}]\le E_{0,\om}[D^+_{T_k}]\le k+1
\end{equation}\abel{harry}
 by Lemma \ref{beatus}.
Therefore, (\ref{snape}) implies
\begin{equation}\label{hagrid}
\lim_{k\to\infty}\frac{1}{k}E_{0,\om}[D^+_{T_k}] = 1.
\end{equation}\abel{hagrid}
By a calculation similar to the one in (\ref{nonumber}), $\wE_0[D_\infty^k]=\wE_0[D_\infty^0]$.
Consequently, we can proceed like in the proof of \cite[Theorem 12]{ze} as follows and get
\begin{eqnarray*}
\wE_0[D_\infty^0]&=& \frac{1}{K}\sum_{k=0}^{K-1}
\wE_0[D_\infty^k]
\ \ge\  \wE_0\left[ \frac{1}{K}\sum_{k=0}^{K-1}D_{T_K}^k\right]\ =\ 
\wEE\left[\frac{1}{K}E_{0,\om}[D^+_{T_K}]\right].
\end{eqnarray*}
Dominated convergence for $K\to\infty$, justified by (\ref{harry}), and (\ref{hagrid}) then yield 
\begin{equation}\label{lilly}
1\le \wE_0[D_\infty^0]\le \wEE\bigg[\sum_{x\in \SL_0}\delta^x\bigg]= L\wEE[\delta^0].
\end{equation}\abel{lilly}
Now consider the event
$S:=\left\{\sum_{x\in \SL_0}\delta^x>\om(0,e_1,1)-\om(0,-e_1,1)\right\}$
that not all the drift contained in the slab $\SL_0$ is stored in the
first cookie at 0. Observe that $\wPP[S]>0$. Indeed, for $L\ge 2$ this follows
from independence of the environment at different sites and for $L=1$ the opposite
would imply $\wEE[\delta^0]\le 1-\kappa$, contradicting (\ref{lilly}).

Now according to (\ref{unopoo})
we have $\wPP$-a.s.\ $P_{0,\om}[A_\ell]=1$. Therefore, by Lemma \ref{D}, $\wPP$-a.s.\
$P_{0,\om}[A_\ell\cap B_\ell]>0$.
Hence, since $\wPP[S]>0$, as shown above,
\begin{eqnarray*}
0&<&\wEE[P_{0,\om}[B_\ell], S]\le \wP_0[D_\infty^0=\om(0,e_1,1)-\om(0,-e_1,1), S]\\
&\le&  \wP_0\bigg[D_\infty^0<\sum_{x\in \SL_0}\delta^x\bigg].
\end{eqnarray*}
Since $D_\infty^0\le\sum_{x\in \SL_0}\delta^x$ anyway, this implies
$\wE_0[D_\infty^0]<L\wEE[\delta^0]=L\EE[\delta^0].$
Along with (\ref{lilly}) this contradicts the assumption $L\EE[\delta^0]\le 1$.
\end{proof}

We conclude with some remarks, discussing the assumption of uniform ellipticity and  some relation
to branching processes with immigration.

\begin{rem}{\rm 1. The following example shows that the
 assumption of uniform ellipticity in Theorems \ref{dd}  and \ref{kellerer} with $Y\ne\Z$ is essential.
Let  $\ell=e_1$, and let $(\om(x))_x$ be i.i.d.\ under $\PP$ with $\PP[\om(0)=\om_+]=1/2=
\PP[\om(0)=\om_-]$, where $\om_+$ and $\om_-$ are such that for all $i\ge 1$, $(\om_\pm(e,i))_{|e|=1}\in (0,1)^{2d}$ is
a probability transition vector with
$\om_\pm(\pm e_2,i)=1-2^{-i}$ and $\om_\pm(e_1,i)\ge \om_\pm(-e_1,i)$.  
Then, by the  Borel-Cantelli lemma,  the walk  will $P_0$-a.s.\ eventually become periodic 
and get stuck on two random sites
$x$ and $x+e_2$ with $\om(x)=\om_+$ and $\om(x+e_2)=\om_-$. Hence it will not be 
transient and might not 
be recurrent to its starting point.
%

2. It is well-known that  recurrence of the simple 
symmetric random walk $(Y_n)_n$ on $\Z$ corresponds to extinction 
of the Galton-Watson process $(Z_m)_m$ with
 geometric(1/2) offspring distribution. Indeed, let the walk start at $Y_0=1$, set $Z_0=1$ and
denote by $Z_m$, $m\ge 1$, the number of transitions
of $(Y_n)_n$ from $m$ to $m+1$ before the walk hits 0.
Since for ERW transitions to the right are more likely than for $(Y_n)_n$,  ERW can be viewed 
as a Galton-Watson process with immigration. 
Pakes \cite[Theorem 1]{Pa71} and 
Zubkov \cite[Theorem 3]{Zu72} showed that adding to each non-empty generation of a 
critical Galton-Watson process 
an i.i.d.\  number of immigrants
 makes it supercritical if the mean number of immigrants is 
above a certain critical threshold. 
This is reminiscent of Theorem \ref{d1}.
However, since the immigration component of the Galton-Watson process
derived from ERW is not independent these results do not directly translate into results for ERW.
}
\end{rem}
\bibliographystyle{amsalpha}

{\sc\small 
Mathematisches Institut\\
Universit\"at T\"ubingen\\
Auf der Morgenstelle 10\\
72076 T\"ubingen, Germany\\
E-Mail: {\rm martin.zerner@uni-tuebingen.de} }
\end{document}